\documentclass[11pt,a4paper]{amsart}

\usepackage{color}

\begin{document}
\newtheorem{theor}{Th\'eor\`eme}
\newtheorem{prop}[theor]{Proposition}
\newtheorem{cor}[theor]{Corollaire}
\newtheorem{lemma}[theor]{Lemme}
\newtheorem{defin}[theor]{D\'efinition}
\newtheorem{conj}{Conjecture}
\newtheorem{rem}[theor]{Remarques}
\newtheorem{nota}[theor]{Notations}

\def\vs{\vskip 4mm}
\def\eq{equation}
\def\cd{\cdot}
\def\cds{\cdots}
\def\ri{\rightline{$\Box$}}
\def\rid{\rightline{$\Box$}}
\def\R{{\mathbb{R} } }
\def\Z{{\mathbb{Z} } }
\def\Q{{\mathbb{Q} } }
\def\P{{\mathbb{P} } }
\def\N{{\mathbb{N} } }
\def\O{\mathcal{O}}
\def\D{\mathcal{D}}
\def\A{\mathcal{A}}
\def\E{\mathcal{E}}
\def\C{{\mathbb{C}}}
\def\widebar{\overline}

\title[Sur la conjecture de Zagier pour $n=4$. II]{Sur la conjecture de Zagier pour $n=4$. II}
\author{Nicusor Dan}
\thanks{Travail r\'ealis\'e avec le support du CNCSIS-UEFISCU, par le
contrat de recherche PN-II-ID-2228/2008.}
\address{Institutul de Matematica al Academiei Romane\\ Calea Grivitei 21 \\Bucuresti 010702, Romania}
\email{ndan@dnt.ro}

\begin{abstract}
We express a general multiple polylogarithm of weight $n$ as an explicit linear combination of
multiple polylogarithms of weight $n$ in $n-2$ variables.
We express a general multiple polylogarithm of weight $4$ as an explicit linear combination of
multiple polylogarithms of type $(3,1)$.
We deduce a $4$ parameters functional equation expressing a certain linear combination of
multiple polylogarithms of type $(3,1)$ as a linear combination of polylogarithms of weight $4$.

\vs

R\'ESUM\'E. On exprime un polylogarithme multiple de poids $n$ g\'en\'eral comme combinaison
lin\'eaire explicite de polylogarithmes multiples de poids $n$ en $n-2$ variables.
On exprime un polylogarithme multiple de poids $4$ g\'en\'eral comme combinaison
lin\'eaire explicite de polylogarithmes multiples de type $(3,1)$.
On d\'eduit une \'equation fonctionelle \`a $4$ param\`etres qui exprime une certaine combinaison lin\'eaire des polylogarithmes multiples de type $(3,1)$ comme combinaison lin\'eaire des polylogarithmes de poids $4$.

\end{abstract}

\maketitle

\section{Introduction}
\label{100}

\vs Soit $n$ un entier positif. Si $a$ et $b$ sont deux points sur une vari\'et\'e complexe $X$ et
$\omega_1, \cdots, \omega_n$ sont des formes holomorphes de d\'egr\'e $1$ sur
$X$, l'int\'egrale it\'er\'ee se d\'efinit par induction sur $n$
par la formule
$$\int_a^b \omega_1 \circ  \cdots \circ  \omega_n = \int_a^b (\int_a^t \omega_1 \circ
\cdots\circ  \omega_{n-1})\omega_n(t).$$

Les polylogarithmes multiples sont d\'efinis comme les int\'egrales it\'er\'ees
$$H(a_0|a_1,\cdots,a_n|a_{n+1})=\int_{a_0}^{a_{n+1}} \frac{dt}{t-a_1} \circ  \frac{dt}{t-a_2} \circ  \cdots \circ
\frac{dt}{t-a_n}.$$

C'est une fonction complexe multivalu\'ee sur l'ensemble des
$(n+2)-$ulpes complexes $(a_0,\cdots,a_{n+1})$ v\'erifiant $a_0
\neq a_1, a_n \neq a_{n+1}$ (pour que l'int\'egrale converge). Elle est invariante par la transform\'ee affine $(a_i)_i \to (\alpha a_i + \beta)_i$, pour des nombres complexes $\alpha \neq 0$ et $\beta$.

Le polylogarithme classique de poids $n$ est la fonction complexe multivalu\'ee sur $\C \setminus \{ 0,1\} $ qui s'\'ecrit $Li_n(z)= \sum_{k=1}^{\infty}\frac{z^k}{k^n}$ sur le disque unit\'e $|z|\leq 1$. On peut prouver que $Li_n(z)=-H(0|1,0, \cdots,0|z)$, donc le polylogarithme classique est le polylogarithme multiple r\'eduit \`a une seule variable.

Soit $E$ un corps. On va d\'efinir les polylogarithmes multiples "a valeurs on $E$". On note $E_{\times}^{n+2}$ l'ensemble des
$(n+2)-$uples $(a_0,\cdots,a_{n+1})$ de $E$ satisfaisant $a_0\neq
a_1$ et $a_n\neq a_{n+1}$. On note $E_{\times}^{n+2} / (E^{\times} \times E)$ le quotient de
$E_{\times}^{n+2}$ par les transformations affines $(a_i)_i \to (\alpha a_i + \beta)_i$. On note $\mathcal{A}_n(E)$ l'espace vectoriel sur $\Q$ ayant comme base les symboles
$[a_0|a_1, \cdots, a_n|a_{n+1}]$ pour $(a_0, \cdots, a_{n+1})\in
E_{\times}^{n+2} / (E^{\times} \times E)$. On d\'efinit ${\mathcal A}_0(E)=\Q$. L'espace vectoriel gradu\'e $\mathcal{A}(E)=\oplus_{n\geq 0} \mathcal{A}_n(E)$ admet une structure de bialg\`ebre. La multiplication est donn\'ee par la formule
\begin{\eq}
\label{900}
[a_0|a_1,\cdots,a_k|a_{k+l+1}]\cdot [a_0|a_{k+1},\cdots,a_{k+l}|a_{k+l+1}]=\sum_{\sigma}[a_0|a_{\sigma(1)},\cdots,a_{\sigma(k+l)}|a_{k+l+1}],
\end{\eq}
o\`u $\sigma$ parcourt les permutations de  l'ensemble $\{ 1,\cdots,k+l \}$ qui pr\'eservent l'ordre dans les sous-ensembles $\{ 1,\cdots,k \}$ et $\{ k+1,\cdots,k+l \}$. La comultiplication est donn\'ee par une formule plus compliqu\'ee ([3]).

On consid\`ere la coalg\`ebre de Lie  $\mathcal{B}(E)=\mathcal{A}(E)/\mathcal{A}_{>0}(E)\cdot\mathcal{A}_{>0}(E)$ des \'el\'ements primitifs. On note $\delta=\oplus_n \delta_n: \mathcal{B} (E)\to \mathcal{B} (E)\otimes \mathcal{B}(E)$ la cod\'erivation. Suivant Zagier, on d\'efinit par induction l'espace vectoriel
$\mathcal{R}_n(E)\subset \mathcal{B}_n(E)$ des
"r\'elations entre polylogarithmes multiples on $E$" et on pose $\mathcal{H}_n(E):=\mathcal{B}_n(E)/\mathcal{R}_n(E)$. On note
toujours $[a_0|a_1, \cdots, a_n|a_{n+1}]$ la classe de
l'\'el\'ement $[a_0|a_1, \cdots, a_n|a_{n+1}]$ modulo $\mathcal{R}_n(E)$.  Le sous-espace
vectoriel $\mathcal{R}_1(E)$ est par d\'efinition engendr\'e
par $[a|z|b] + [b|z|c] = [a|z|c]$ pour  $z, a, b, c \in E$ v\'erifiant $z\neq a,b,c$. On calcule
$\mathcal{H}_1(E)=E^{\times}\otimes_{\Z}{\Q}$. On note $\mathcal{K}_n(E)$ le noyau de l'aplication $(pr \otimes pr)\circ \delta_n : \mathcal{B}_n(E) \to (\mathcal{H}(E) \otimes \mathcal{H}(E))_n$, o\`u $pr: \mathcal{B}_k(E)\to \mathcal{H}_k(E)$ a \'et\'e d\'ej\`a d\'efini pour $k<n$. Soit $t$ une variable. On d\'efinit $\mathcal{R}_n(E)$ comme le sous-espace vectoriel engendr\'e par $\alpha(1) -
\alpha(0)$ pour tous les \'el\'ements $\alpha$ de $\mathcal{K}_n(E(t))$ pour lesquels $\alpha(1)$ et $
\alpha(0)$ sont bien d\'efinies. L'application  $(pr \otimes pr)\circ \delta_n$ se factorise \`a une application $\delta_n : \mathcal{H}_n(E) \to (\mathcal{H}(E) \otimes \mathcal{H}(E))_n$ qui fait de $\mathcal{H}$ une coalg\`ebre de Lie gradu\'ee.

La conjecture de Zagier ([4]) affirme que la valeur en $s=n$ de la fonction z\^eta de Dedekind d'un corps de nombre est le d\'eterminant d'une matrice dont les termes sont des polylogarithmes de poids $n$ \'evalu\'ees dans des \'el\'ements du corps en question. Apr\`es des travaux de Goncharov et Zagier, la conjecture se r\'eduit \`a une conjecture qui affirme que le r\'egulateur de Beilinson est combinaison lin\'eaire des polylogarithmes. On peut prouver que le r\'egulateur de Beilinson est combinaison lin\'eaire des polylogarithmes multiples. Il reste \`a trouver des formules exprimant les polylogarithmes multiples (en $n$ variables) comme combinaisons lin\'eaires de polylogarithmes (polylogarithmes multiples en $1$ variable). L'article [1] donne une pr\'esentation synthetique de ces reductions. Dans le pr\'esent article on parcourt les premiers deux pas de la strategie: passer de $n$ variables \`a $n-2$ variables.

\section{Le th\'eor\`eme pour $n$ g\'en\'eral}

Soit $n$ un entier positif. On introduit une g\'en\'eralisation l\'eg\`ere de la notion de polylogarithme multiple. Soient $a_0, a_1, \cdots, a_{n+1}, x$ des \'el\'ements de ${\P}^1(\C)$ v\'erifiant $a_0\neq a_1, a_0\neq x, a_n\neq a_{n+1}, x\neq a_{n+1}$. On choisit $\omega(a_i,x)$ l'unique forme diff\'erentielle de d\'egr\'e 1 holomorphe sur ${\P}^1(\C)-\{ a_i, x\}$ qui est nulle si $a_i=x$ et qui a un p\^ole d'ordre $1$ de r\'esidu $+1$ en $a_i$ et un p\^ole d'ordre $1$ de r\'esidu $-1$ en $x$ si $a_i \neq x$. On d\'efinit
$$H(a_0|a_1,\cdots,a_n//x|a_{n+1})=\int_{a_0}^{a_{n+1}} \omega(a_1,x) \circ  \cdots \circ  \omega(a_n,x)$$
Est une fonction invariante pour l'action de $PGL(2, \C)$ sur $((a_i)_i,x)$. Le passage entre cette fonction et la fonction ant\'erieure est clair:
\begin{\eq}
\label{1000}
H(a_0|a_1,\cdots,a_n|a_{n+1})= H(a_0|a_1,\cdots,a_n//\infty|a_{n+1}),
\end{\eq}
$$ H(a_0|a_1,\cdots,a_n//x|a_{n+1})=H((a_0-x)^{-1}|(a_1-x)^{-1},\cdots,(a_n-x)^{-1}|(a_{n+1}-x)^{-1})$$
si $x\neq \infty$. Soit $y$ un autre \'el\'ement de ${\P}^1(\C)$. En \'ecrivant $\omega(a_i,x) = \omega(a_i, y)-\omega(x,y)$ pour tout $0\leq i \leq n$ et en d\'ev\'elopant de mani\`ere multilin\'eaire on peut \'ecrire le polylogarithme multiple $H(a_0|a_1,\cdots,a_n//x|a_{n+1})$ comme somme altern\'ee des polylogarithmes multiples $H(a_0|\cdots//y|a_{n+1})$. On va exploiter cette r\'elation.

Soit $E$ un corps. On transpose les consid\'erations ci-dessus aux \'el\'ements de $\mathcal{H}(E)$. Soient $a_0, a_1, \cdots, a_{n+1}, x$ des \'el\'ements distincts de ${\P}^1(E)$. En analogie avec $(\ref{1000})$, on d\'efinit l'\'el\'ement $[a_0|a_1,\cdots,a_n//x|a_{n+1}]$ dans $\mathcal{H}_n(E)$ comme \'etant $[a_0|a_1,\cdots,a_n|a_{n+1}]$ si $x=\infty$ et $[(a_0-x)^{-1}|(a_1-x)^{-1},\cdots,(a_n-x)^{-1}|(a_{n+1}-x)^{-1}]$ si $x\neq \infty$. Pour $1\leq i \leq n$ et pour $I$ un sous-ensemble de l'ensemble $\{1, \cdots, n\}$ cont\'enant $i$, on d\'efinit
$A([a_0|a_1,\cdots,a_n//x|a_{n+1}],i,I)$ comme le symbole $[a_0|a_1,\cdots,a_n//x|a_{n+1}]$ dans lequel on remplace $a_j$ par $a_i$ dans toutes les positions $j \in I$. On d\'efinit
$$B([a_0|a_1,\cdots,a_n//x|a_{n+1}],i)=\sum_{I} (-1)^{|I|} A([a_0|a_1,\cdots,a_n//x|a_{n+1}],i,I),$$
la somme \'etant prise sur toutes les sous-ensembles $I$ de l'ensemble $\{1, \cdots, n\}$ cont\'enant $i$ et ayant le cardinal $|I|\geq 2$. Les consid\'erations du paragraphe pr\'ec\'edent apliqu\'ees \`a $y=a_i$ sug\'erent la r\'elation dans $\mathcal{H}_n(E)$:
\begin{\eq}
\label{2000}
[a_0|a_1,\cdots,a_n//x|a_{n+1}]+[a_0|a_1,\cdots,a_{i-1},x,a_{i+1},\cdots, a_n//a_i|a_{n+1}]
\end{\eq}
$$=B([a_0|a_1,\cdots,a_n//x|a_{n+1}],i)$$

C'est une r\'elation dans $\mathcal{H}_n(E)$ qu'on v\'erifie facilement par induction sur $n$.

On consid\`ere un entier $0\leq s \leq n$. Par induction sur $s$, en utilisant des r\'elations $(\ref{900})$, on peut prouver facilement
$$[a_0|y,\cdots,y, b_{s+1}, b_{s+2}, \cdots, b_n//x|a_{n+1}]=(-1)^s\sum_J C_J,$$
o\`u $J$ parcourt les sous-ensembles de cardinal $s$ de l'ensemble $\{2, \cdots, n\}$ et $C_J$ d\'esigne le symbole $[a_0|b_{s+1}, \cdots//x|a_{n+1}]$ dans lequel les positions $J$ sont ocup\'ees par les $y$ et les positions restantes par $b_{s+2}, \cdots, b_n$ dans cette ordre. On aplique cela \`a $A([a_0|a_1,\cdots,a_n//x|a_{n+1}],i,I)$, \`a $y=a_i$ et \`a $s$ le plus petit entier pour lequel $\{1, \cdots, s\} \subset I$. On obtient une pr\'esentation de $A([a_0|a_1,\cdots,a_n//x|a_{n+1}],i,I)$ comme somme altern\'ee explicite de $[a_0|b_1,\cdots,b_n//x|a_{n+1}]$ o\`u les $(b_j)_j$ sont parmi les $(a_j)_j$, il y a au moins deux $a_i$ parmi les $(b_j)_j$ et aucun d'entre eux sur la premi\`ere position. Il est facile de prouver que, dans ces notations
$$[a_0|b_1,\cdots,b_n//x|a_{n+1}]=[a_i|b_1,\cdots,b_n//x|a_{n+1}]-[a_i|b_1,\cdots,b_n//x|a_0].$$
On peut donc \'ecrire $A([a_0|a_1,\cdots,a_n//x|a_{n+1}],i,I)$ comme somme altern\'ee explicite des polylogarithmes multiples en $\leq n-2$ variables. On remplace dans $(\ref{2000})$ et on obtient
$$[a_0|a_1,\cdots,a_n//x|a_{n+1}]+[a_0|a_1,\cdots,a_{i-1},x,a_{i+1},\cdots, a_n//a_i|a_{n+1}]$$
$$=D([a_0|a_1,\cdots,a_n//x|a_{n+1}],i),$$
o\`u $D([a_0|a_1,\cdots,a_n//x|a_{n+1}],i)$ est une somme altern\'ee explicite des polylogarithmes multiples en $\leq n-2$ variables.

Soient $1 \leq i<j \leq n$ deux entiers. En appliquant trois fois la r\'elation pr\'ec\'edente pour les substitutions $x\to a_i$, $a_i \to a_j$, $a_j \to x$, on obtient
\begin{\eq}
\label{3000}
[a_0|a_1,\cdots,a_n//x|a_{n+1}]+[a_0|a_1,\cdots,a_{i-1},a_j,a_{i+1},\cdots,a_{j-1},a_i,a_{j+1}, \cdots, a_n//x|a_{n+1}]
\end{\eq}
$$=D([a_0|a_1,\cdots,a_n//x|a_{n+1}],i) - D([a_0|a_1,\cdots,a_{i-1},x,a_{i+1},\cdots, a_n//a_i|a_{n+1}],j)$$
$$+D([a_0|a_1,\cdots,a_{i-1},x,a_{i+1},\cdots, a_{j-1},a_i, a_{j+1}, \cdots, a_n//a_j|a_{n+1}],i).$$

On d\'eduit que, pour toute permutation $\sigma$ de l'ensemble $\{ 1,\cdots, n\}$,
$$[a_0|a_{\sigma(1)},\cdots,a_{\sigma(n)}//x|a_{n+1}] - {\mathrm sign}(\sigma)[a_0|a_1,\cdots,a_n//x|a_{n+1}]$$ est une somme altern\'ee explicite des polylogarithmes multiples en $\leq n-2$ variables.

On suppose maintenant $n\geq 3$. La r\'elation $(\ref{900})$ pour $k=2$, $l=n-2$ dans $\mathcal{H}_n(E)$ donne
\begin{\eq}
\label{4000}
\sum_{\sigma\in S}[a_0|a_{\sigma(1)},\cdots,a_{\sigma(n)}|a_{n+1}]=0,
\end{\eq}
o\`u $S$ est l'ensemble des permutations de  l'ensemble $\{ 1,\cdots,n \}$ qui pr\'eservent l'ordre dans les sous-ensembles $\{ 1,2 \}$ et $\{ 3,\cdots,n \}$. Comme $\sum_{\sigma \in S} {\mathrm sign}(\sigma) = [n/2] \neq 0$, o\`u $[x]$ d\'enote la partie enti\`ere du nombre r\'eel $x$, on peut trouver une combinaison lin\'eaire des r\'elations $(\ref{3000})$ qui, adition\'ee \`a la r\'elation $(\ref{4000})$, exprime $[n/2][a_0|a_1,\cdots,a_n//x|a_{n+1}]$ comme combinaison lin\'eaire des polylogarithmes multiples en $\leq n-2$ variables. On peut faire cela de maniere explicite. On note, pour $1\leq i<j\leq n$, $A_{i,j}=[a_0|a_{\sigma(1)},\cdots,a_{\sigma(n)}|a_{n+1}]$ pour la permutation $\sigma$ qui met $1$ sur la position $i$ et $2$ sur la position $j$. On note $R(i-1,j|i,j)$ la r\'elation $(\ref{3000})$ pour $A_{i-1,j}+A_{i,j}$ et de m\^eme pour $R(i,j|i,j+1)$. On doit consid\'erer la somme de la r\'elation $(\ref{4000})$ avec la combinaison lin\'eaire $\sum_{1<i<j}c(i-1,j|i,j)R(i-1,j|i,j)+\sum_{1<j<n}c(1,j|1,j+1)R(1,j|1,j+1)$, o\`u $c(i-1,j|i,j)=-1$ si $j-i$ est impair et $0$ sinon et o\`u $c(1,j|1,j+1)=(-1)^j([n/2]-[j/2])$. On d\'eduit

\begin{theor}
\label{theor1}
Soit $n\geq 3$ un entier. Soit $E$ un corps. Soient $a_0, a_1, \cdots, a_{n+1}$ des \'el\'ements distincts de $E$. Alors $[a_0|a_1, \cdots, a_n|a_{n+1}]$ est combinaison lin\'eaire explicite des polylogarithmes multiples de poids $n$ en $\leq n-2$ variables dans $\mathcal{H}_n(E)$.
\end{theor}

\par{\bf Remarques: }
1) Francis Brown m'a comuniqu\'e d'avoir prouv\'e que $[a_0|a_1, \cdots, a_n|a_{n+1}]$ peut s'\'ecrire comme combinaison lin\'eaire des polylogarithmes multiples de poids $n$ en $\leq n-2$ variables dans $\mathcal{H}_n(E)$. Il n'est pas clair si sa m\'ethode peut \^etre rendue explicite.

2) Si $n$ est impair on peut obtenir une formule plus simple. On utilise la r\'elation $(\ref{900})$ pour $k=1$, $l=n-1$ pour obtenir une r\'elation $(\ref{4000})$ pour $S$ l'ensemble des permutations de  l'ensemble $\{ 1,\cdots,n \}$ qui pr\'eservent l'ordre dans le sous-ensemble $\{ 2,\cdots,n \}$. On a $\sum_{\sigma \in S} {\mathrm sign}(\sigma) = 1$. On note, pour $1<i<n$, $R_i$ la r\'elation $(\ref{3000})$ pour $[a_0|a_2,\cdots, a_i, a_1, a_{i+1},\cds, a_n|a_{n+1}]+[a_0|a_2,\cdots, a_{i+1}, a_1, a_{i+2},\cds, a_n|a_{n+1}]$. On doit consid\'erer la r\'elation $(\ref{4000}) - R_2 -R_4 - \cds - R_{n-1}$.

\section{Le cas $n=4$}
Quand $n=4$, les polylogarithmes multiples de poids $4$ en $2$ variables sont $[x,y]_{3,1}=[0|x,0,0,y|1]$, $[x,y]_{2,2}=[0|x,0,y,0|1]$, $[x,y]_{1,3}=[0|x,y,0,0|1]$ et les polylogarithmes multiples de poids $4$ en $1$ variable sont les polylogarithmes $[x]_4=[0|x,0,0,0|1]$. Il existent des formules explicites exprimant les fonctions $[x,y]_{2,2}$ et $[x,y]_{1,3}$ comme combinaisons lin\'eaires des fonctions $[x,y]_{3,1}$ et polylogarithmes. Si on remplace ces formules dans le Th\'eor\`eme \ref{theor1} on obtient:
\begin{theor}
\label{theor2}
Soit $E$ un corps. Soient $a,b,c,d,e,f$ des \'el\'ements distincts de $E$. On note pour simplifier $[\cdot, \cdot] = [\cdot, \cdot]_{3,1}$, $abc = \frac{a-c}{b-c}$, $abcd=\frac{(a-c)(b-d)}{(a-d)(b-c)}$ et de m\^eme pour les autres combinaisons. Alors on a dans $\mathcal{H}_4(E)$:
$$[a|b,c,d,e|f]=\phi(a,b,c,d,e)-\phi(f,b,c,d,e),$$
o\`u
$$2\phi(a,b,c,d,e)=[aed,ced]-[ecd,acd]-2[cad,ead]+2[acd,bcd]-[bad,cad]-[cbd,abd]$$
$$-[cab,eab]+[aeb,ceb]-[dab,cab]+[acb,dcb]-[dab,eab]+[aeb,deb]-[ace,dce]+[dae,cae]$$
$$-[cbe,abe]-[bae,cae]+[bde,ade]+[abe,dbe]-[dca,dcae]-[dcbe,dcba]-[dcab,dca]$$
$$+[dcb,dcba]+[cda,cdb]+[ecd,ecda]-[cea,ced]-[ecad,eca]+[ecab,eca]-[ecb,ecba]$$
$$-[cea,ceb]+[ecda,ecdb]+[dbac,dbae]+[dbea,dbe]-[dba,dbae]-[bde,bda]+[bdc,bda]$$
$$-[dba,dbac]-[dbca,dbc]+[ebca,ebc]-[bec,bea]-[eba,ebac]+[ebdc,ebda]+[bea,bed]$$
$$+\gamma(a,b,c,d,e),$$
et $\gamma(a,b,c,d,e)$ est une combinaison lin\'eaire explicite des polylogarithmes $[x]_4$, o\`u chaque $x$ est fraction rationelle en $a,b,c,d,e$.
\end{theor}

\par{\bf Remarques: }
1) Dans l'article [1], on a obtenu une autre pr\'esentation de $[a|b,c,d,e|f]$ comme combinaison lin\'eaire explicite de fonctions $[x,y]_{3,1}$ et $[z]_4$. La comparaison des deux formules donne une \'equation fonctionelle \`a $4$ param\`etres qui exprime une certaine combinaison lin\'eaire des fonctions de type $[x,y]_{3,1}$ comme combinaison lin\'eaire des polylogarithmes $[z]_4$.

2) Dans l'article [1], on a montre comment la conjecture de Zagier pour $n=4$ se r\'eduit \`a une conjecture de Goncharov ([2]), qui pr\'edit  qu'une certaine combinaison lin\'eaire des fonctions de type $[x,y]_{3,1}$ \`a $3$ param\`etres peut \^etre \'ecrite comme combinaison lin\'eaire des polylogarithmes $[z]_4$. Il serait int\'eressant de comparer la combinaison lin\'eaire des fonctions de type $[x,y]_{3,1}$ de la remarque pr\'ec\'edente avec celle de la conjecture de Goncharov. Si la deuxi\`emme peut \^etre d\'eduite de la premi\`ere, on aurait prouv\'e la conjecture de Zagier pour $n=4$.

3) Dans la formule de [1], on a \'ecrit $[a|b,c,d,e|f]=F(a,b,c,d,e)-F(f,b,c,d,e)$. En plus, la fonction $F(a,b,c,d,e)$ est invariante \`a la permutation cyclique $a\to b\to c\to d\to e\to a$. Il n'est pas clair si $F(a,b,c,d,e)$ est \'egal \`a $\phi(a,b,c,d,e)$. Si cela est vrai, l'\'equation fonctionelle de la remarque 1) est somme de deux \'equations fonctionelles \`a $3$ param\`etres.

\vs \par{\bf Bibliographie:}

[1]: N. Dan: Sur la conjecture de Zagier pour $n=4$, arXiv:0809.3984 [math.KT] (2008)

[2]: A. B. Goncharov: Polylogarithms and motivic Galois group,
Proc. Sympos. Pure Math., vol. 55, Part 2, AMS, Providence, RI
(1994), p. 43-96

[3]: A. B. Goncharov: Multiple polylogarithms and mixed Tate motives, arXiv:math/0103059 [math.AG] (2001)

[4]: D. Zagier: Polylogarithms, Dedekind zeta functions and the
algebraic K-theory of fields, Progr. Math, vol. 89(1991), p.
391-430

\end{document}